\documentclass[11pt]{article}
\usepackage[pagebackref=true,colorlinks,linkcolor=green,citecolor=magenta]{hyperref}
\usepackage{amssymb,amsmath}
\usepackage{amsmath}
\usepackage{dsfont}
\usepackage{amsfonts}
\newtheorem{precor}{{\bf Corollary}}

\newtheorem{precon}{{\bf Conjecture}}

\newtheorem{prealphcon}{{\bf Conjecture}}

\newtheorem{predefin}{{\bf Definition}}

\newtheorem{preexm}{{\bf Example}}

\newtheorem{preappl}{{\bf Application}}

\newtheorem{prelem}{{\bf Lemma}}

\newtheorem{preproof}{{\bf Proof.\ }}

\newtheorem{prethm}{{\bf Theorem}}

\newtheorem{prealphthm}{{\bf Theorem}}

\newenvironment{alphthm}{\begin{prealphthm}{\hspace{-0.5
               em}{\bf.\ }}}{\end{prealphthm}}

\newtheorem{prealphlem}{{\bf Lemma}}

\newtheorem{prepro}{{\bf Proposition}}

\newtheorem{preprb}{{\bf Problem}}

\newtheorem{prerem}{{\bf Remark}}

\newtheorem{preapp}{{\bf Application}}

\newtheorem{prequ}{{\bf Question}}

\newtheorem{preclaim}{{\bf Claim}}

\def\conct[#1,#2]{\mbox {${#1} \leftrightarrow {#2}$}}
\def\dconct[#1,#2]{\mbox {${#1} \rightarrow {#2}$}}
\def\deg[#1,#2]{\mbox {$d_{_{#1}}(#2)$}}
\def\mindeg[#1]{\mbox {$\delta_{_{#1}}$}}
\def\maxdeg[#1]{\mbox {$\Delta_{_{#1}}$}}
\def\outdeg[#1,#2]{\mbox {$d_{_{#1}}^{^+}(#2)$}}
\def\minoutdeg[#1]{\mbox {$\delta_{_{#1}}^{^+}$}}
\def\maxoutdeg[#1]{\mbox {$\Delta_{_{#1}}^{^+}$}}
\def\indeg[#1,#2]{\mbox {$d_{_{#1}}^{^-}(#2)$}}
\def\minindeg[#1]{\mbox {$\delta_{_{#1}}^{^-}$}}
\def\maxindeg[#1]{\mbox {$\Delta_{_{#1}}^{^-}$}}

\def\dre[#1,#2,#3]{\mbox {${\cal E}^{^{#3}}(#1,#2)$}}
\def\var[#1,#2]{\mbox {${\rm Var}_{_{#1}}(#2)$}}
\def\ls[#1]{\mbox {$\xi^{^{#1}}$}}
\def\hom[#1,#2]{\mbox {${\rm Hom}({#1},{#2})$}}
\def\onvhom[#1,#2]{\mbox {${\rm Hom^{v}}(#1,#2)$}}
\def\onehom[#1,#2]{\mbox {${\rm Hom^{e}}(#1,#2)$}}
\def\core[#1]{\mbox {$#1^{^{\bullet}}$}}
\def\cay[#1,#2]{\mbox {${\rm Cay}({#1},{#2})$}}
\def\sch[#1,#2,#3]{\mbox {${\rm Sch}({#1},{#2},{#3})$}}
\def\cays[#1,#2]{\mbox {${\rm Cay_{s}}({#1},{#2})$}}
\def\dirc[#1]{\mbox {$\stackrel{\rightarrow}{C}_{_{#1}}$}}
\def\cycl[#1]{\mbox {${\bf Z}_{_{#1}}$}}

\begin{document}
\footnotetext[1]{The research of Hossein Hajiabolhassan is supported by ERC advanced grant GRACOL.}
\begin{center}
{\Large \bf  A Note on Altermatic Number}\\
\vspace{0.3 cm}
{\bf Meysam Alishahi$^\dag$ and Hossein Hajiabolhassan$^\ast$\\
{\it $^\dag$ School of Mathematical Sciences}\\
{\it University of Shahrood, Shahrood, Iran}\\
{\tt meysam\_alishahi@shahroodut.ac.ir}\\
{\it $^\ast$ Department of Mathematical Sciences}\\
{\it Shahid Beheshti University, G.C.}\\
{\it P.O. Box {\rm 19839-69411}, Tehran, Iran}\\
{\tt hhaji@sbu.ac.ir}\\
{\it $^\ast$ School of Mathematics}\\
{\it Institute for Research in Fundamental Sciences~(IPM)}\\
{\it P.O. Box {\rm 19395-5746}, Tehran, Iran}\\
}
\end{center}
\begin{abstract}
\noindent
In view of Tucker's lemma (an equivalent combinatorial version of the Borsuk-Ulam theorem), the present authors (2013) introduced
the $k^{th}$ altermatic number of a graph $G$ as a tight lower bound for the chromatic number of $G$.
In this note, we present a purely combinatorial proof for this result.\\

\end{abstract}
\section{Introduction}
Throughout the paper, the set $\{1,2,\ldots, n\}$ is denoted by $[n]$.
A {\it hypergraph} ${\cal H}$ is an ordered pair $(V({\cal H}),E({\cal H}))$,
where $V({\cal H})$ and $E({\cal H})$ are the vertex set and the hyperedge set of ${\cal H}$, respectively.
Hereafter, all hypergraphs are  simple, i.e.,
$E({\cal H})$ is a family of distinct nonempty subsets of $V({\cal H})$.
A proper $t$-coloring of a hypergraph ${\cal H}$, is a mapping $c:V({\cal H})\longrightarrow [t]$
such that no hyperedge is monochromatic, i.e.,
$|c(e)|>1$ for any hyperedge $e\in E({\cal H})$.
The minimum integer $t$ such that ${\cal H}$ admits a $t$-coloring
is called the {\it chromatic number} of ${\cal H}$ and is denoted by $\chi({\cal H})$.
For a hypergraph containing some hyperedge of cardinality $1$, we define its chromatic number to be infinite.

For a  vector $X=(x_1,x_2,\ldots,x_n)\in \{R,0,B\}^n$,
a subsequence $x_{a_1}, x_{a_2},\ldots,x_{a_t}$ ($1\leq {a_1}<{a_2}<\cdots<{a_t}\leq n$) of nonzero terms of $X$
is called an {\it alternating subsequence} of $X$ if any
two consecutive terms in this subsequence are different. We denote by $alt(X)$ the length of a longest alternating subsequence of $X$.
Moreover, we define $alt(0,0,\ldots,0)=0$. Also,
we denote the number of nonzero terms of $X$ by $|X|$.
For instance, if
$X=(R,R,B,B,0,R,0,R,B)$, then $alt(X)=4$ and $|X|=7$.
For an $X=(x_1,x_2,\ldots,x_n)\in\{R,0,B\}^n$, define
$X^R=\{i\ :\ x_i=R\}$ and $X^B=\{i\ :\ x_i=B\}$.
Note that if we consider a vector $X=(x_1,x_2,\ldots,x_n)$, then one can obtain $X^R$ and $X^B$, and conversely.
Therefore, by abuse of notation, we can set $X=(X^R,X^B)$.
Throughout the paper, we use interchangeably these representations, i.e.,
$X=(x_1,x_2,\ldots,x_n)$ or $X=(X^R,X^B)$.
For $X=(X^R,X^B), Y=(Y^R,Y^B)\in \{R, 0, B\}^n$, we write
$X\subseteq Y$, if $X^R \subseteq Y^R$ and $X^B \subseteq Y^B$.
Note that if $X\subseteq Y$, then every alternating subsequence of $X$
is an alternating subsequence of $Y$, and subsequently, $alt(X)\leq alt(Y)$.
Also, if the first nonzero term of $X$ is $R$ (resp. $B$),
then every alternating subsequence of $X$ with the maximum length begins with $R$ (resp. $B$), and moreover, we can conclude that $X^R$ (resp. $X^B$) contains the smallest integer of $X^R\cup X^B$.

Let $L_{V({\cal H})}=\{v_{i_1}<v_{i_2}\ldots<v_{i_n}\ :\ (i_1,i_2,\ldots,i_n)\in S_n\}$ be the set of all linear orderings of
the vertex set of hypergraph ${\cal H}$, where $V({\cal H})=\{v_1,v_2,\ldots,v_n\}$.
For any $X=(x_1,x_2,\ldots,x_n)\in \{R, 0, B\}^n$ and any linear ordering
$\sigma: v_{i_1}<v_{i_2}\ldots<v_{i_n}\in L_{V({\cal H})}$, define $X^R_\sigma=\{v_{i_j}\ :\ x_j=R\}$,  $X^B_\sigma=\{v_{i_k}\ :\ x_k=B\}$, and $X_\sigma=(X^R_\sigma,X^B_\sigma)$.
Note that for $V=[n]$ and $I:1<2<\cdots<n$, we have $X^R=X^R_I$, $X^B=X^B_I$, and $X=X_I$.
Also, set ${\cal H}_{|_{X_\sigma}}$ to be the hypergraph with the vertex set $X^R_\sigma\cup X^B_\sigma$  and the edge set
$$E({\cal H}_{|_{X_\sigma}})=\left\{A\in E({\cal H}):\  \ A\subseteq X^R_\sigma\ {\rm or}\ A\subseteq X^B_\sigma\right\}.$$

For a hypergraph ${\cal H}=(V({\cal H}),E({\cal H}))$, the {\it general Kneser graph} ${\rm KG}({\cal H})$ has all hyperedges of
${\cal H}$ as vertex set and two vertices of ${\rm KG}({\cal H})$ are adjacent if the corresponding hyperedges are disjoint.
The hypergraph ${\cal H}$ provides a {\it Kneser representation} for a graph $G$ whenever $G$ and ${\rm KG}({\cal H})$ are isomorphic. It is simple to see that a graph $G$ has various Kneser representations.
For any $\sigma\in L_{V({\cal H})}$ and positive integer $k$,
define $alt_\sigma({\cal H}, k)$
to be the largest integer $t$ such that there exists an
$X\in\{R,0,B\}^n$ with  $alt(X)=t$
and that the chromatic number of  ${\rm KG}({\cal H}_{|X_\sigma})$ is at most $k-1$ (For $k=1$, it means that ${\cal H}_{|X_\sigma}$ contains no
hyperedge). Now define $$alt({\cal H},k)=\min\left\{alt_\sigma({\cal H},k):\ \sigma\in L_{V({\cal H})}\right\}.$$

Let $G$ be a graph and $k$ be a positive
integer such that $1\leq k\leq\chi(G)+1$.
The {\it $k^{th}$ altermatic number} of $G$, $\zeta(G,k)$, is defined as follows
$$
\zeta(G,k)=\displaystyle\max_{{\cal H}}
\left\{|V({\cal H})|-alt({\cal H},k)+k-1: {\rm KG}({\cal H})\longleftrightarrow G\right\}
$$
where ${\rm KG}({\cal H})\longleftrightarrow G$ means there are some homomorphisms from $G$ to $H$ and also from $H$ to $G$.

One can see that for $k= \chi({\rm KG}({\cal H}))+1$, we have $alt({\cal H},k)=|V|$; and consequently, $\chi({\rm KG}({\cal H}))=\zeta(G,k)$.
In~\cite{2013arXiv1302.5394A}, in view of Tucker's lemma, it was shown that the $k^{th}$ altermatic number of a graph is a tight lower bound for its chromatic number.
\begin{alphthm}\label{combin}{\rm \cite{2013arXiv1302.5394A}}
For any graph $G$ and positive integer $k$, where $k\leq\chi(G)+1$, we have
$\chi(G)\geq \zeta(G,k).$
\end{alphthm}
Furthermore, it was shown~\cite{2013arXiv1302.5394A} that the first altermatic number can be considered as an improvement
of the Dol'nikov-K{\v{r}}{\'{\i}}{\v{z}}'s lower bound~\cite{MR953021, MR1081939} for the chromatic number of
general Kneser hypergraphs.  Also, we should mention that by an improvement of Gale's lemma,  in~\cite{2014arXiv1403.4404A}, it was shown  that the
first altermatic number provides tight lower bound for some well-known topological parameter which is related to graphs via the Borsuk-Ulam theorem.

For more details about Tucker's lemma and Gale's lemma, we refer readers to~\cite{2014arXiv1403.4404A,MR1988723}. In~\cite{2013arXiv1306.1112M}, Meunier showed that
it is a hard problem to determine the exact value of the first altermatic number of a graph. Precisely, he proved  that for a hypergraph ${\cal H}$ and
a permutation $\sigma$ of $V({\cal H})$,
it is an NP-hard problem to specify $alt_\sigma({\cal H},1)$.

\section{A Combinatorial Proof of Theorem~\ref{combin}}
For two positive integers $m$ and $n$ where $m\geq 2n$,
the usual Kneser graph ${\rm KG}(m,n)$ is a graph whose vertex set consists of all
$n$-subsets of $[m]$ and two vertices are adjacent if the corresponding subsets are disjoint.
Kneser  in 1955~\cite{MR0068536} conjectured that $\chi({\rm KG}(m,n))\geq m-2n+2$.
In 1978~\cite{MR514625}, Lov\'asz by using algebraic topology, proved Kneser conjecture. Next Schrijver~\cite{MR512648} introduced  the Schrijver graph ${\rm SG}(m,n)$ as a subgraph  of ${\rm KG}(m,n)$ and proved that it is  critical  and has the same chromatic number as ${\rm KG}(m,n)$.
It is proved in~\cite{2013arXiv1302.5394A} that Theorem~\ref{combin} (for $k=1,2$ ) implies these results. In what follows, we present the proof of Theorem~\ref{combin}.
It should be mentioned that the proof relies on an interesting idea used by Matou{\v{s}}ek~\cite{MR2057690}
to present a combinatorial proof of Lov\'asz-Kneser theorem.


\noindent{\bf Proof of Theorem~\ref{combin}.}
On the contrary, suppose $\zeta(G,k) > \chi(G)$. Consider a hypergraph
${\cal H}$ such that ${\rm KG}({\cal H})$ is isomorphic to $G$ and that
$$\zeta(G,k)\geq  |V({\cal H})|-alt({\cal H},k)+k-1=|V({\cal H})|-alt_\sigma({\cal H},k)+k-1> \chi(G)$$
where $\sigma\in L_{V({\cal H})}$.
Without loss of generality and for the simplicity of notations, we may assume that
$V=[n]$ and $\sigma=I:1<2<\cdots<n$.
Let  $h: V({\rm KG}({\cal H}))=E({\cal H}) \longrightarrow \{1,2,\ldots, n-alt({\cal H},k)+k-2\}$ be a proper coloring of
${\rm KG}({\cal H})$. For any subset $M\subseteq V({\cal H})$, we define $\bar{h}(M)=\max\{h(A):\ A\subseteq M,\ A\in E({\cal H})\}$. If there is no $A\subseteq M$, where $A\in E({\cal H})$,
then set $\bar{h}(M)=0$.
For $X=(X^R,X^B)\in \{R,0,B\}^n$, set $\bar{h}(X)=\max\{\bar{h}(X^R),\bar{h}(X^B)\}$.
Define a map $\lambda:\ \{R,0,B\}^n \longrightarrow \{\pm1,\pm2,\ldots,\pm n\}$ as follows
\begin{itemize}
\item If $X=(X^R,X^B)\in \{R,0,B\}^n$ and
      $alt(X)\leq alt_{I}({\cal H},k)$, set
$$
\lambda(X)=\left\{
\begin{array}{cl}
+alt(X)+1 & {\rm if}\ X^B=\varnothing \ {\rm or}\ \min (X^R \cup X^B)\in X^R\\
-alt(X)-1   & {\rm otherwise}
\end{array}\right.
$$
\item If $X=(X^R,X^B)\in \{R,0,B\}^n$ and
      $alt(X)\geq alt_{I}({\cal H},k)+1$, set
$$
\lambda(X)=\left\{
\begin{array}{cl}
alt_{I}({\cal H},k)+\bar{h}(X)-k+2 & {\rm if}\ \bar{h}(X)=\bar{h}(X^R)\\
-(alt_{I}({\cal H},k)+\bar{h}(X)-k+2) & {\rm if}\  \bar{h}(X)=\bar{h}(X^B)
\end{array}\right.
$$
\end{itemize}
Since $h$ is a proper coloring, one can see that the map $\lambda$ is well-defined. Moreover,
it is straightforward to check that for any two ordered pairs $(A_i, B_i)\subseteq (A_j, B_j)$, we have
$\lambda((A_i, B_i))+\lambda((A_j, B_j))\not =0$. Also,
in view of the definition of $alt_{I}({\cal H},k)$,  if $alt(X)\geq alt_{I}({\cal H},k)+1$, then
the chromatic number of ${\rm KG}(F_{|_{X}})$ is at least $k$, and consequently, $|\lambda(X)|\geq alt_{I}({\cal H},k)+2$.
Also, by the definition of $\lambda$, one can see that $\lambda((\varnothing,\varnothing))=1$.
In the sequel, we show that there exists a graph $H$ with a unique vertex of
degree one and any other vertex of degree $2$, which is impossible. This contradicts our assumption that $h$ is a proper coloring.
For any subset $A\subseteq [n]$, define $-A=\{-t: t\in A\}$. Note that if $A=\varnothing$,
then $-A=\varnothing$.
A permissible sequence is a sequence
$(A_0, B_0)\subseteq (A_1, B_1)\subseteq \cdots \subseteq(A_m, B_m)$ of disjoint ordered pairs of $\{R,0,B\}^n$
such that
for any $0\leq i \leq m$, we have $|A_i|+|B_i|=i$ and that $A_m\cup -B_m\subseteq \{\lambda((A_0, B_0)),\lambda((A_1, B_1)),\ldots, \lambda((A_m, B_m))\}$. By definition, we have  $A_0,=B_0=\varnothing$ and $0\leq m\leq n$. Also, one can see that $(\varnothing,\varnothing)$ is a permissible sequence.

Set the vertex set of the graph $H$ to be the set of all permissible sequences. Now we introduce
the edge set of $H$. For the permissible sequence $(\varnothing,\varnothing)$, we define its unique neighbor to be
$(\varnothing,\varnothing)\subseteq (\{1\},\varnothing)$, which is a permissible sequence.
We assign two neighbors to any other permissible sequence as follows. Moreover, we show that this assignment is symmetric,
i.e., $H$ is an undirected graph.
Consider a permissible sequence
$(A_0, B_0)\subseteq (A_1, B_1)\subseteq \cdots \subseteq(A_m, B_m)$  ($m\geq 1$) and set $\lambda_i=\lambda((A_i, B_i))$ for any $0\leq i \leq m$.
By the definition of $\lambda$, it is clear that if $(A,B)\subseteq (A',B')$,
then  $|\lambda(A,B)|\leq |\lambda(A',B')|$. Therefore,
in view of the definition of permissible sequence,
one of the following conditions holds
\begin{enumerate}
\item[(i)]  There exists a unique  integer $0\leq i < m$ such that $\lambda_i=\lambda_{i+1}$.
\item[(ii)] There exists a unique integer $0\leq i \leq m$ such that $\lambda_i\not \in A_m\cup -B_m$.
\end{enumerate}

Note that in case~(i),  $i =0$ is~not possible since for any $j>0$, $|\lambda_j|>1$.
If $\lambda_i=\lambda_{i+1}$ for $1\leq i < m$, then we define two neighbors of  $(A_0, B_0)\subseteq (A_1, B_1)\subseteq \cdots \subseteq(A_m, B_m)$ as follows.
\begin{enumerate}
\item $(A'_0, B'_0)\subseteq (A'_1, B'_1)\subseteq \cdots \subseteq(A'_m, B'_m)$, where
$(A'_i, B'_i)=(A_{i-1}\cup(A_{i+1}\setminus A_i), B_{i-1}\cup(B_{i+1}\setminus B_i))$
and for any $r\not = i$, $(A'_r, B'_r)=(A_r, B_r)$.

\item If $i<m-1$, then define the other neighbor to be
$(A''_0, B''_0)\subseteq (A''_1, B''_1)\subseteq \cdots \subseteq(A''_m, B''_m)$, where $(A''_{i+1}, B''_{i+1})=(A_{i}\cup(A_{i+2}\setminus A_{i+1}), B_{i}\cup(B_{i+2}\setminus B_{i+1}))$
and for any $r\not = i+1$, $(A''_r, B''_r)=(A_r, B_r)$. Otherwise, if $i=m-1$, define the other neighbor to be
$(A_0, B_0)\subseteq (A_1, B_1)\subseteq \cdots \subseteq(A_{m-1}, B_{m-1})$.
\end{enumerate}
One can check that both of neighbors are permissible. Now suppose that
there exists an integer $0\leq i \leq m$ such that $\lambda_i\not \in A_m\cup -B_m$. Define the neighbors as follows
\begin{enumerate}
\item $(A'_0, B'_0)\subseteq (A'_1, B'_1)\subseteq \cdots \subseteq(A'_m, B'_m)\subseteq(A'_{m+1}, B'_{m+1})$, where
for any $0\leq r\leq m$, $(A'_r, B'_r)=(A_r, B_r)$ and if $\lambda_i>0$, then
$(A'_{m+1},B'_{m+1})=(A_{m}\cup\{\lambda_i\},B_{m})$, otherwise, $(A'_{m+1},B'_{m+1})=(A_{m},B_{m}\cup\{-\lambda_i\})$.

\item If $1\leq i\leq m-1$, then define the other neighbor to be
$(A''_0, B''_0)\subseteq (A''_1, B''_1)\subseteq \cdots \subseteq(A''_m, B''_m)$, where $(A''_i, B''_i)=(A_{i-1}\cup(A_{i+1}\setminus A_i), B_{i-1}\cup(B_{i+1}\setminus B_i))$
and for any $r\not = i$, $(A''_r, B''_r)=(A_r, B_r)$. If $i=m$, define the second neighbor to be
$(A_0, B_0)\subseteq (A_1, B_1)\subseteq \cdots \subseteq(A_{m-1}, B_{m-1})$. Otherwise, for $i=0$, consider
$(B_0, A_0)\subseteq (B_1, A_1)\subseteq \cdots \subseteq(B_m, A_m)$ as the second neighbor.
\end{enumerate}
Note that all neighbors are permissible.
Also, one can check that the aforementioned assignment is symmetric, which completes the proof.
\hfill$\blacksquare$\\

\noindent {\bf Acknowledgement:}
This paper was written while Hossein Hajiabolhassan was visiting School of Mathematics, Institute for Research in Fundamental Sciences~(IPM).
He acknowledges the support of IPM (No. $94050128$).

\def\cprime{$'$} \def\cprime{$'$}


\begin{thebibliography}{10}

\bibitem{2013arXiv1302.5394A}
M.~Alishahi and H.~Hajiabolhassan.
\newblock On the chromatic number of general kneser hypergraphs.
\newblock {\em Journal of Combinatorial Theory, Series B}, 115:186 -- 209,
  2015.

\bibitem{2014arXiv1403.4404A}
M.~{Alishahi} and H.~{Hajiabolhassan}.
\newblock {Hedetniemi's Conjecture Via Altermatic Number}.
\newblock {\em ArXiv e-prints}, March 2014.

\bibitem{MR953021}
V.~L. Dol{\cprime}nikov.
\newblock A combinatorial inequality.
\newblock {\em Sibirsk. Mat. Zh.}, 29(3):53--58, 219, 1988.

\bibitem{MR0068536}
M.~Kneser.
\newblock Ein {S}atz \"uber abelsche {G}ruppen mit {A}nwendungen auf die
  {G}eometrie der {Z}ahlen.
\newblock {\em Math. Z.}, 61:429--434, 1955.

\bibitem{MR1081939}
I.~K{\v{r}}{\'{\i}}{\v{z}}.
\newblock Equivariant cohomology and lower bounds for chromatic numbers.
\newblock {\em Trans. Amer. Math. Soc.}, 333(2):567--577, 1992.

\bibitem{MR514625}
L.~Lov{\'a}sz.
\newblock Kneser's conjecture, chromatic number, and homotopy.
\newblock {\em J. Combin. Theory Ser. A}, 25(3):319--324, 1978.

\bibitem{MR1988723}
J.~Matou{\v{s}}ek.
\newblock {\em Using the {B}orsuk-{U}lam theorem}.
\newblock Universitext. Springer-Verlag, Berlin, 2003.
\newblock Lectures on topological methods in combinatorics and geometry,
  Written in cooperation with Anders Bj{\"o}rner and G{\"u}nter M. Ziegler.

\bibitem{MR2057690}
J.~Matou{\v{s}}ek.
\newblock A combinatorial proof of {K}neser's conjecture.
\newblock {\em Combinatorica}, 24(1):163--170, 2004.

\bibitem{2013arXiv1306.1112M}
F.~{Meunier}.
\newblock {Colorful Subhypergraphs in Kneser Hypergraphs}.
\newblock {\em Electron. J. Combin.}, 21(1):\ Research Paper \#P1.8, 13 pp.
  (electronic), 2014.

\bibitem{MR512648}
A.~Schrijver.
\newblock Vertex-critical subgraphs of {K}neser graphs.
\newblock {\em Nieuw Arch. Wisk. (3)}, 26(3):454--461, 1978.

\end{thebibliography}
\end{document}